\documentclass[twocolumn]{autart2}
\usepackage[utf8]{inputenc}
\usepackage[english]{babel}
\usepackage{graphicx}      
\usepackage{cite}
\usepackage{kantlipsum,cuted}

\usepackage{graphics} 
\usepackage{epsfig} 
\usepackage{graphicx,color}
\usepackage{subfigure}
\usepackage{ushort}
\usepackage{caption}
\usepackage[overload]{empheq}

\usepackage{multirow}
\usepackage{amssymb,amsfonts,enumerate}
\usepackage[subnum]{cases}
\newcommand{\R}{{\mathbb R}}
\newtheorem {remark}{Remark}
\newtheorem {defi}{Definition}

\newcommand{\mysgn}{\mathrm{sgn^+}}
\usepackage{pgf,tikz}
\usetikzlibrary{arrows}
\definecolor{zzttqq}{rgb}{0.6,0.2,0}
\definecolor{aqaqaq}{rgb}{0.63,0.63,0.63}
\definecolor{qqzzcc}{rgb}{0,0.6,0.8}
\definecolor{dcrutc}{rgb}{0.86,0.08,0.24}
\definecolor{cqcqcq}{rgb}{0.75,0.75,0.75}
\definecolor{qqwuqq}{rgb}{0,0.39,0}
\definecolor{qqzzff}{rgb}{0,0.6,1}
\definecolor{xdxdff}{rgb}{0.49,0.49,1}
\definecolor{uuuuuu}{rgb}{0.27,0.27,0.27}
\definecolor{qqccqq}{rgb}{0,0.8,0}
\definecolor{ttffqq}{rgb}{0,.5,.5}
\definecolor{ffttww}{rgb}{0.8,0.8,0.8}
\definecolor{qqqqff}{rgb}{0,0,1}
\definecolor{rosso}{rgb}{1,0,0}

\begin{document}
\begin{frontmatter}
\title{Estimation and control of oscillators through short-range noisy proximity measurements}
\author[Naples]{Francesco Lo Iudice}
\ead{francesco.loiudice2@unina.it}    
\author[Sevilla]{Jos\'e \'Angel Acosta}
\ead{jaar@us.es}              
\author[Naples]{Franco Garofalo}
\ead{franco.garofalo@unina.it}
\author[Naples]{Pietro DeLellis} 
\ead{pietro.delellis@unina.it}
\address[Naples]{Department of Electrical Engineering and Information Technology,\\ University of Naples Federico II, Naples 80125, Italy}
\address[Sevilla]{Department of Automatic Control and System Engineering, University of Sevilla, Sevilla, 41004, Spain }
\begin{abstract}
In this paper, we present a novel estimation and control strategy to balance a formation of discrete-time oscillators on a circle. We consider the case in which each oscillator only gathers noisy proximity measurements, whose range is lower than the desired spacing along the circle, implying total disconnectedness of the balanced formation. These restrictions pose relevant challenges that are overcome through the symbiotic combination of an estimator that borrows tools from interval analysis and a three-level bang-bang controller. We prove that the formation can be balanced, with an accuracy that can be regulated by tuning a controller parameter. The effectiveness of the proposed strategy is further illustrated through a set of numerical simulations.
\end{abstract}
\end{frontmatter}
\section{Introduction}
\vspace{-0.3cm}
Coordinating the motion of multi-agent systems is a relevant issue in very diverse fields of science and engineering spanning from biology to robotics \cite{WaLa:91,ReBe:05,MaCo:07,abpo17}. In formation control, most works rely on the agents being able to directly measure their relative position \cite{KnCh:16,CaYu:13,BaAr:97}. However, when only distance measurements are available, coordination becomes significantly harder\cite{cayu11,jide17,yean17} as the intrinsic ambiguity of these measurement calls for complementing the controller with an estimator able to reconstruct the agents' relative position. While distance measurements can be obtained with sensors based on different technology, a common trait among these is having limited range \cite{KiKi:18,BeBl:02}. Accounting for the sensors' range through proximity communication rules \cite{Wa:06,SaKo:08} is necessary when budget constraints do not allow the deployment of long range sensors, and poses additional challenges to the estimation and control strategy \cite{DeGa:15,dega15}.
\vspace{-0.3cm}

Achieving a balanced circular formation has emerged as a paradigmatic formation control problem \cite{ChZh:11,MaBr:04,WaXi:13,SoLi:16,YuXi:18}. It has been tackled by assuming that the relative position is measurable \cite{MaBr:04,SmBr:05,KiSu:07,ChZh:13} and under rather strong connectivity assumptions. An all-to-all connectivity was assumed in \cite{SePa:08,sewu14}, while fixed and connected graphs were considered in \cite{jide17,yean17,MaBr:04,momi09}. However, in the presence of proximity rules, none of the above results can be applied, as the relative positions are not available, and the measurement flow is intermittent. Recently, a discontinuous control law was proposed to solve this control problem \cite{ChZh:11,ChZh:13} assuming joint connectivity of the proximity graph. However, when the sensor range is too short compared to the desired distance along the circle, this assumption becomes too restrictive.
\vspace{-0.3cm}

In this paper, we devise an estimation and control strategy capable of balancing the formation without making any connectivity assumption. The agents are first order integrators on the unit circle, and the problem is directly stated in discrete-time in view of the implementation.
Our strategy determines a relative motion between a randomly elected pacemaker and the other agents, thus triggering a cascade in which each agent $i$ identifies its closest follower through an estimation algorithm and then varies its speed through a bang-bang controller to adjust its distance with respect to its follower. In turn, this speed variation induces a relative motion between $i$ and the next agent thus allowing the latter to identify $i$ as its follower. The cascade only stops when the formation is achieved. By tuning the controller parameters, it is possible to regulate the pace of the multi-agent system, the balancing accuracy, and the convergence speed towards the steady-state formation. Summing up, differently from the existing literature,
\vspace{-0.3cm}

\noindent (1) no assumption on the connectivity of the graph is required. Convergence is proved assuming that the detecting distance is lower than the desired spacing, thus implying a disconnected steady-state topology;
\newline (2) neither the absolute or the relative position among the agents is measured: our strategy only requires a (noisy) proximity measurement;
\newline (3) we provide bounds on the convergence time and on the accuracy of the formation balancing as explicit functions of the controller parameters, that can be then regulated depending on the performance required.

\section{Problem statement}\label{sec:statement}
\subsection{Mathematical preliminaries and notation.}
Given an interval $J \subset \mathbb{R}$, we denote its infimum $\ushort J$, its supremum $\bar J$ and its width $w(J):=\bar{J}-\ushort{J} \ \in \mathbb{R}$. The Minkowski sum between two intervals $X,Y\in\R$ is $\lbrace x + y\ |\ x \in X, y \in Y \rbrace$. As a scalar can be seen as a degenerate interval, all sums in this paper are to be intended as Minkowski sums. Given $\lambda$ intervals $J_1,\ldots,J_\lambda$, the infimum and the supremum of the interval hull $H=\mathrm{hull}_l\left\{J_l\right\}$ are given by $\ushort{H}={\mathrm{inf}}_l \lbrace \ushort{J}_{l}\rbrace$ and 
$\bar{H}={\mathrm{sup}}_l \lbrace \bar{J}_{l}\rbrace$, respectively.
We define the function
\begin{equation}
\mysgn(x):=
\begin{cases}
1 \ \mathrm{if} \ x\geq0,\\
0 \ \mathrm{otherwise}.
\end{cases}
\end{equation}
The floor $\lfloor x \rfloor$ and ceiling $\lceil x \rceil$ functions, associate to each $x\in\mathbb{R}$ the largest integer not greater than $x$ and the smallest integer not less than $x$, respectively. Given $x\in\mathbb{R}$, we denote by $\mathrm{rem}(x)$ the unique solution for $r$ to the equation
$x=2\pi q + r,$ where $-\pi\le r< \pi$, $q\in\mathbb{Z}$.
\subsection{Agent dynamics and control goal}
We consider $N$ oscillators on a circle whose angular position dynamics are described by
\begin{equation}
\theta_i(k+1) = \theta_i(k)+\omega+u_i(k), \ \forall i=1,\dots,N\label{eq:multi_agent}
\end{equation}
where $\omega$ is the natural angular speed, and $u_i(k)$ is the control input at time $k$. Introducing the relative angular position $\theta_{ij}(k) := \theta_i(k)-\theta_j(k)$, we have
\begin{equation}\label{eq:rel_ang_pos}
\theta_{ij}(k+1) = \theta_{ij}(k)+u_{ij}(k),  
\end{equation}
where $u_{ij}(k) := u_i(k)-u_j(k)$.
Without loss of generality, we assume that $\theta_{ij}(0) \in [-\pi,\pi)$ for all $i,j=1,...,N$. 
Also, let the relative phase difference be defined as $\vartheta_{ij}(k):=\mathrm{rem}(\theta_{ij}(k)) \in [-\pi, \pi)$.
\begin{defi}
Given a scalar $\varepsilon>0$, we say that the multi-agent system \eqref{eq:multi_agent} achieves a $\varepsilon$-partially balanced circular formation if, for all $\theta_{ij}(0)$, $i,j=1,\ldots,N,i\ne j$, 
\begin{equation}\label{eq:eps_balanced}
\limsup_{k\rightarrow\infty}\left|\vartheta_{ij}(k)-\psi\right|\le \varepsilon,
\end{equation}
for all $(i,j)\in\left\{(1,2),\ldots,(N-1,N),(N,1)\right\}$, and where 
$\psi:=2\pi/N$ is the spacing distance.
\end{defi}
We aim at designing a control strategy $u_i$, $i=1,\ldots,N$, such that \eqref{eq:eps_balanced} holds for some finite $\varepsilon$ and whose parameters can be tuned to make this $\varepsilon$ smaller. We will assume that that the agents can only rely on intermittent, short-range and noisy proximity measurements. In particular,
\begin{enumerate}[(a)]
\item  we measure the angular distance $\alpha_{ij}(k)$ instead of $\vartheta_{ij}(k)$, defined as
$\alpha_{ij}(k) := |\vartheta_{ij}(k)|$;

\item for each pair of agents, the measurement $y_{ij}(k)$ of $\alpha_{ij}(k)$ is only available if $\alpha_{ij}(k)\leq \theta_{\max}>0$;
\item the measurement $y_{ij}(k)$, when available, is affected by a bounded noise $\nu_{ij}(k)$;
\item the detecting distance $\theta_{\max}$ is lower than the desired spacing distance $\psi$.
\end{enumerate}

This setting forces each agent to estimate the relative angular position with respect to the others before deciding the control input and implies that the output of system \eqref{eq:rel_ang_pos} be
\begin{equation}\label{eq:measurement_equation} 
y_{ij} (k) = 
\begin{cases}
\alpha_{ij}(k)+\nu_{ij}(k) & \mbox{if} \ \alpha_{ij}(k) \in I,\\ 
\mbox{no measure} & \mbox{otherwise}, 
\end{cases}
\end{equation}
where $I:= [0,\theta_{\max}]$ and $\nu_{ij}(k)$ is the measurement noise whose amplitude is bounded by $\varphi$, for all $i,j=1,\ldots,N$. Notice that, even in absence of noise, $\vartheta_{ij}(k)$ could not be directly computed from $y_{ij}(k)$, as two phase differences with opposite signs are compatible with the same measurement $y_{ij}(k)$. This intermittent measurement flow can be described through the time-varying graph $\mathcal{G}(k)=\{\mathcal{V},\mathcal{E}(k)\}$, where $\mathcal{V}=\{1,\ldots,N \}$ and $(i,j)\in\mathcal{E}(k)$ if $\alpha_{ij}(k)\in I$.
Therefore, point (d) implies that, when the desired spacing $\psi$ is achieved, the proximity graph is not connected, that is, $\mathcal{E}(k)=\emptyset$. Accordingly, in our estimation and control design we cannot rely on connectivity. 
\vspace{-0.2cm}
\section{Strategy for estimation and control}\label{sec:est_contr}
\vspace{-0.2cm}
Our strategy for achieving a partially balanced circular formation requires labeling each agent and randomly electing a {\em pacemaker}, from now on denoted as agent 1, whose motion will not be affected by that of its peers. The remainder of the agents implement an estimation procedure based on that presented in \cite{DeGa:15} that combines the information brought by the measurements with that brought by the knowledge of the dynamics to build a finite multi-interval set, $\Gamma^{ij}(k|k)$, where the relative phase $\vartheta_{ij}(k)$ among the agents falls. This estimate is leveraged by the agents to identify their closest follower $i - 1$, defined as $ i-1 := {\mbox{argmin}}_j \lbrace \vartheta_{ij} \geq 0 \rbrace$ and then exploited by a decentralized bang-bang control law that achieves a balanced circular formation by allowing each agent to be \emph{pushed} by its closest follower.
\vspace{-0.3cm}
\subsection{Preliminaries}
\vspace{-0.3cm}
We start our preliminary considerations by exploiting the information that each measure $y_{ij}(k)$ brings on the angular distance $\alpha_{ij}(k)$. For all $k$, we know that
\vspace{-0.1cm}
\begin{equation}
\alpha_{ij}(k) \in 
\begin{cases}
\Upsilon_{ij}(k)  \   \mbox{if a measure is available,}\label{eq:def_Upsilon} \\
I^c  \ \mathrm{otherwise},
\end{cases}
\end{equation}
where $I^c := (\theta_{\max},\pi]$ and $\Upsilon_{ij}(k) := [\max \lbrace y_{ij}(k)-\varphi,0 \rbrace,\min \lbrace y_{ij}(k)+\varphi,\theta_{\max} \rbrace]\subseteq I.$
As $\vartheta_{ij}(k)$ is related to $\alpha_{ij}(k)$ through the absolute value function, at each time instant $k$ relation \eqref{eq:def_Upsilon} allows to identify two intervals in which $\vartheta_{ij}(k)$ falls. By considering the information brought by the knowledge of the dynamics of agents $i$ and $j$, our estimation strategy reduces these two intervals to one, recursively shrinks its width, and extracts a scalar estimate $\hat \vartheta_{ij}(k)$ of $\vartheta_{ij}(k)$. This estimate is then exploited to achieve our control goal.
Hence, at each time instant $k$, our knowledge on $\vartheta_{ij}(k)$ will be represented by the set $\Gamma^{ij}(k|k)$ which, in general, is composed of the union of two intervals $\Gamma_1^{ij}(k|k) \subset [0,\pi)$ and $\Gamma_2^{ij}(k|k)\subset[-\pi, 0]$. In what follows, we denote its hull by $H^{ij}(k|k)$, which represents an overestimate of the uncertainty on $\vartheta^{ij}(k)$. The following definition is introduced 1) to provide the conditions guaranteeing that a generic agent $i$ has identified its follower, and 2) to introduce the notation $k_i$ for the first time-instant in which agent $i$ has identified its follower. 
\vspace{-0.5cm}
\begin{defi}\label{def:follower}
Agent $i\ne 1$ identifies its closest follower at time $k_i$ if $k_i$ is the smallest integer ensuring there exists $k\leq k_i$ such that
\vspace{-0.5cm}
\begin{subequations}\label{eq:follower}
\begin{align}[left ={\empheqlbrace}]\label{eq:follower_a}
\ushort H^{i,i-1}(k|k)>0,\\
\label{eq:follower_b}
\bar H^{i,i-1}(k|k) < \ushort{\Gamma}_1^{ij}(k|k), \ \forall j \neq i-1.
\end{align}
\end{subequations}
\end{defi}
The function $\mathcal{I}_i(k)$ tracks which agents, except the pacemaker, have already identified their closest follower, i.e.
\begin{equation}\label{eq:indic_fun}
{\mathcal I}_i(k) = 
\begin{cases}
1 \ & \forall k \geq k_i \\
0 \ & \forall k< k_i,
\end{cases}
\qquad  i=2,\ldots, N .
\end{equation}
\subsection{Decentralized estimation and control laws}\label{subsec:est_contr}
To achieve an $\varepsilon$-partially balanced circular formation we employ the following three-level bang bang controller:
\begin{subequations}\label{eq:contr_strat}
\begin{align}[left ={u_{i}(k) = \empheqlbrace}]
\label{eq:our_control_law_others_a}
\omega_0 + K\mathrm{sgn^+}( \psi- &\hat \vartheta_{i,i-1}(k))  \\
 &\text{if} \ {\mathcal I}_i(k) = 1,i\ne 1\nonumber\\
\label{eq:our_control_law_leader}
\omega_0 \qquad \qquad \qquad \, &\text{if }i=1\\
\label{eq:our_control_law_others_c}
0 \qquad \qquad \qquad  \ \ &\mbox{otherwise}
\end{align}
\end{subequations}
where $\hat \vartheta_{i,i-1}(k)$ is a scalar estimate of $\vartheta_{ij}(k)$ made by agent $i$, and formally defined in eq. \eqref{eq:scal_est} and $K$ is a tunable control parameter. 

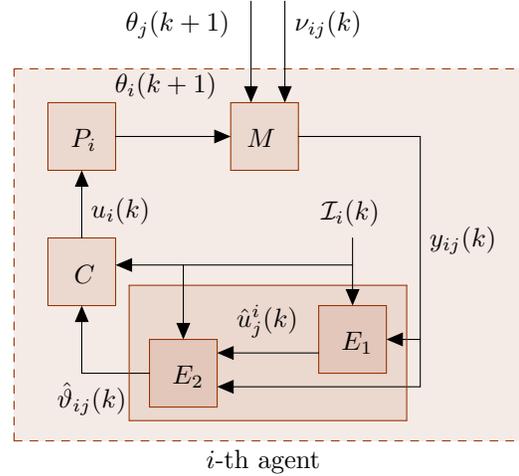
\begin{figure}
\centering
\begin{tikzpicture}[line cap=round,line join=round,>=triangle 45,x=0.9cm,y=0.9cm]
\clip(1,1.4) rectangle (9,8.8);
\fill[color=zzttqq,fill=zzttqq,fill opacity=0.1] (2,6) -- (3,6) -- (3,7) -- (2,7) -- cycle;	
\fill[color=zzttqq,fill=zzttqq,fill opacity=0.1] (4.7,6) -- (5.7,6) -- (5.7,7) -- (4.7,7) -- cycle;
\fill[color=zzttqq,fill=zzttqq,fill opacity=0.1] (2,4) -- (3,4) -- (3,5) -- (2,5) -- cycle;
\fill[color=zzttqq,fill=zzttqq,fill opacity=0.1] (3.5,2.5) -- (4.5,2.5) -- (4.5,3.5) -- (3.5,3.5) -- cycle;
\fill[color=zzttqq,fill=zzttqq,fill opacity=0.1] (6,3) -- (7,3) -- (7,4) -- (6,4) -- cycle;
\fill[color=zzttqq,fill=zzttqq,fill opacity=0.1] (3.2,2.3) -- (7.3,2.3) -- (7.3,4.3) -- (3.2,4.3) -- cycle;
\fill[color=zzttqq,fill=zzttqq,fill opacity=0.1] (1.5,2) -- (1.5,7.5) -- (9,7.5) -- (9,2) -- cycle;
\draw [color=zzttqq,dashed] (1.5,2) -- (1.5,7.5);
\draw [color=zzttqq,dashed] (1.5,2) -- (9,2);
\draw [color=zzttqq,dashed] (9,7.5) -- (9,2);
\draw [color=zzttqq,dashed] (1.5,7.5) -- (9,7.5);
\draw [color=zzttqq] (2,6)-- (3,6);
\draw [color=zzttqq] (3,6)-- (3,7);
\draw [color=zzttqq] (3,7)-- (2,7);
\draw [color=zzttqq] (2,7)-- (2,6);
\draw [color=zzttqq] (4.7,6)-- (5.7,6);
\draw [color=zzttqq] (5.7,6)-- (5.7,7);
\draw [color=zzttqq] (5.7,7)-- (4.7,7);
\draw [color=zzttqq] (4.7,7)-- (4.7,6);
\draw [->] (3,6.5) -- (4.7,6.5);
\draw (2.2,6.75) node[anchor=north west] {$P_i$};
\draw (4.8,6.75) node[anchor=north west] {$M$};
\draw (5.7,6.5) -- (7.5,6.5);
\draw [color=zzttqq] (2,4)-- (3,4);
\draw [color=zzttqq] (3,4)-- (3,5);
\draw [color=zzttqq] (3,5)-- (2,5);
\draw [color=zzttqq] (2,5)-- (2,4);
\draw [->] (2.5,5) -- (2.5,6);
\draw [color=zzttqq] (3.5,2.5)-- (4.5,2.5);
\draw [color=zzttqq] (4.5,2.5)-- (4.5,3.5);
\draw [color=zzttqq] (4.5,3.5)-- (3.5,3.5);
\draw [color=zzttqq] (3.5,3.5)-- (3.5,2.5);
\draw [color=zzttqq] (6,3)-- (7,3);
\draw [color=zzttqq] (7,3)-- (7,4);
\draw [color=zzttqq] (7,4)-- (6,4);
\draw [color=zzttqq] (6,4)-- (6,3);
\draw [->] (7.5,3.5) -- (7,3.5);
\draw (7.5,6.5)-- (7.5,2.8);
\draw [->] (7.5,2.8) -- (4.5,2.8);
\draw [->] (6,3.3) -- (4.5,3.3);
\draw [->] (2.5,3) -- (2.5,4);
\draw (3.5,3)-- (2.5,3);
\draw [color=zzttqq] (3.2,2.3)-- (7.3,2.3);
\draw [color=zzttqq] (7.3,2.3)-- (7.3,4.3);
\draw [color=zzttqq] (7.3,4.3)-- (3.2,4.3);
\draw [color=zzttqq] (3.2,4.3)-- (3.2,2.3);
\draw [->] (5,8.5) -- (5,7);
\draw [->] (5.5,8.5) -- (5.5,7);
\draw (5.51,8.5) node[anchor=north west] {$\nu_{ij}(k)$};
\draw (3.01,8.5) node[anchor=north west] {$\theta_{j}(k+1)$};
\draw (7.5,5.3) node[anchor=north west] {$y_{ij}(k)$};
\draw (2.85,7.58) node[anchor=north west] {$\theta_i(k+1)$};
\draw (2.25,4.73) node[anchor=north west] {$C$};
\draw (3.7,3.25) node[anchor=north west] {$E_2$};
\draw (6.2,3.75) node[anchor=north west] {$E_1$};
\draw (4.65,4.15) node[anchor=north west] {$\hat u^{i}_{j}(k)$};
\draw (2,3) node[anchor=north west] {$\hat{\vartheta}_{ij}(k)$};
\draw (4.2,2) node[anchor=north west] {$i$-th agent};
\draw (5.9,5.7) node[anchor=north west] {$\mathcal{I}_{i}(k)$};
\draw (2.5,5.73) node[anchor=north west] {$u_i(k)$};
\draw [->] (6.5,5) -- (6.5,4);
\draw [->] (6.5,4.6) -- (3,4.6);
\draw [->] (4,4.6) -- (4,3.5);
\end{tikzpicture}
\caption{$P_i$ are the dynamics of the $i$-th agent given in \eqref{eq:multi_agent}; $E_1$ is the estimator employed by $i$ estimate $u_{j}$, $j\ne i$, see \eqref{eq:u_est}; $E_2$ is the estimator of $\theta_{i,i-1}$ given in \eqref{eq:est_strat}, \eqref{eq:scal_est}; $C$ is the bang-bang controller described in \eqref{eq:contr_strat}; $M$ represents the measurement equation given in \eqref{eq:measurement_equation}.}
\label{fig:scheme}
\end{figure}

While in general estimation algorithms rely on knowledge of both the dynamics and the input signals, the hypothesis that agent $i$ is aware of the control signal $u_{j}$ exerted by another agent $j$ is not compatible with the need of deploying a decentralized strategy. Therefore, to obtain $\hat \vartheta_{i,i-1}(k)$, every agent will perform its own interval estimate of $u_{j}$ according to the following rules:
\begin{strip}\label{eq:uij}
\fbox{\begin{minipage}{\textwidth}
\begin{subequations}
\label{eq:u_est}
\begin{align}[left ={\hat u^i_j(k) = \empheqlbrace}]
\omega_0, \qquad \qquad \, \ &\text{if } \mathcal{I}_i(k) = 1, 
(i,j)\notin\mathcal{E}(k), j=i-1\\
[\omega_0, \omega_0+K], \, \ &\text{if } \mathcal{I}_i(k) = 1,
 (i,j)\in\mathcal{E}(k), j=i-1\\
[0, \omega_0 + K],\quad\, &\text{otherwise} 
\end{align}
\end{subequations}
\end{minipage}}
\end{strip}
\noindent From $\hat u^i_j(k)$ agent $i$ derives an estimate of the relative angular velocity $\hat u_{ij}(k):= u_i(k)-\hat u^i_j(k)$ of $u_{ij}(k)$ which is then employed to dynamically propagate the multi-interval $\Gamma^{ij}(0|-1)$ defined in \eqref{eq:c_iniz} according to eq. \eqref{eq:a_priori}. At each time instant $\Gamma^{ij}(k+1|k)$ is then intersected with the multi-interval resulting from the measurement procedure, see \eqref{eq:a_post_nonid_a}-\eqref{eq:a_post_fol}. Equation \eqref{eq:a_post_nonfol} prescribes that, as soon as each agent has identified its follower $i-1$, it ceases to estimate the position of all other agents $j\neq i-1$ as our control law is designed so that each agent is \emph{pushed} by its closest follower. 
\begin{remark}
In \cite{DeGa:15} it is shown that under the same assumptions on $\theta_{ij}(0)$ made in this paper, the multi-interval $\Gamma^{ij}(k|k)$ always contains the true value of $\theta_{ij}(k)$ for all $k$. The only assumption made in \cite{DeGa:15} that is not fulfilled in this paper is the knowledge of $u_{ij}(k)$, of which, in this case, we perform an (interval) estimate $\hat u_{ij}(k)$.
\end{remark}
\begin{strip}
\fbox{\begin{minipage}{\textwidth}
\begin{subequations}
\label{eq:est_strat}
\begin{equation}\label{eq:a_priori}
\Gamma^{ij}(k+1|k) = \Gamma^{ij}(k|k) + \hat u_{ij}(k) \ \forall k \ \mathrm{and} \ \forall j,
\end{equation}
  \begin{align}[left ={\Gamma^{ij}(k|k) = \empheqlbrace}]
    I^c \cup -I^c \quad \qquad \qquad \qquad \qquad \qquad
    &\text{if } \mathcal{I}_i(k) = 0, (i,j)\notin\mathcal{E}(k), \label{eq:a_post_nonid_a}\\
    \Gamma^{ij}(k|k-1) \cap (\Upsilon^{ij}(k)\cup -\Upsilon^{ij}(k)) \ 
    &\text{if }\mathcal{I}_i(k) = 0, (i,j)\in\mathcal{E}(k)\text{ or }
    \mathcal{I}_i(k) = 1, (i,j)\in\mathcal{E}(k), j=i-1\label{eq:a_post_nonid_b}\\
   \Gamma^{ij}(k|k-1) \cap (I^c \cup -I^c), \qquad\quad\ \ \,
   &\text{if } \mathcal{I}_i(k) = 1, (i,j)\notin\mathcal{E}(k), j=i-1 \label{eq:a_post_fol}\\
   \emptyset \ \qquad \qquad \qquad \qquad \qquad \qquad \qquad &\text{if } \mathcal{I}_i(k) = 1,  j \neq i-1,  \qquad \qquad \ \ \label{eq:a_post_nonfol}
   \end{align}
\begin{equation}\label{eq:c_iniz}
\Gamma^{ij}(0|-1) = [-\pi, \ \pi)
\end{equation}  
  \end{subequations}
\end{minipage}}
\end{strip}
The scalar estimate of $\vartheta_{i,i-1}$ needed in eq. \eqref{eq:our_control_law_others_a} is
\begin{equation}\label{eq:scal_est}
\hat{\vartheta}_{i,i-1}(k) =\bar H^{i,i-1}(k|k).
\end{equation}
A concise schematic of our estimation and control strategy is illustrated in Figure \ref{fig:scheme}.

\noindent\textbf{Assumptions:} in proving convergence of our estimation and control strategy, we make use of four assumptions.
\vspace{-0.3cm}
%
\begin{enumerate}
\item
$|\nu_{ij}(k)| \leq \varphi$ for all $k$, with $\varphi\ge 0$ being a known constant;
\item $|\vartheta_{ij}(0)| \in [\min\lbrace4\varphi+2\omega_0 +2 K,\theta_{\max}\rbrace,\pi], \  \forall i,j = 1,...,N$, $i\neq j$;
\item $2\omega_0+2K<\theta_{\max}$;
\item $\omega_0>0$ and $K>0$.
\end{enumerate}
\begin{remark}
Note that Assumptions 2, 3, and 4 depend on the parameters $\omega_0$ and $K$ of the controller, which can be therefore employed to enforce their fulfillment. Namely,
Assumption 2 implies that, at time $k=0$, the agents must be sufficiently separated to allow an estimate to be recovered before overtaking may occur. Assumption 3 implies that the sampling time must be sufficiently small if compared to the maximum possible agents' relative speed.
\end{remark}
To facilitate the reading of all the following lemmas and theorems, all the symbols contained in their statements are summarized in Table I.
\begin{lem}\label{lem:th_in_H}
Let Assumptions 1-4 hold. If $u_{ij}(k)\in\hat u_{ij}(k)$ for all $k\ge0$, then $\theta_{ij}(k)\in \Gamma^{ij}(k|k)$, for all $k \ge 0$.
\end{lem}
\vspace{-0.7cm}
\begin{pf}
Assume that at $k-1$ we have $\theta_{ij}(k-1) \in \Gamma^{ij}(k-1|k-1)$. Let us define the multi-interval
$G(k|k-1):=\Gamma^{ij}(k-1|k-1) + u_{ij}(k-1).$
From \cite{DeGa:15}, we have $\theta_{ij}(k)\in G(k|k-1)$. From the hypothesis, we have that $u_{ij}(k)\in\hat u_{ij}(k) \ \forall k$, and then $G(k|k-1)\subseteq \Gamma^{ij}(k|k-1)$. Computing $G(k|k)$ with the laws that update $\Gamma^{ij}(k|k-1)$ to $\Gamma^{ij}(k|k)$, see equations  \eqref{eq:a_post_nonid_a}-\eqref{eq:a_post_fol}, and from the properties of interval intersection, we get that $G(k|k) \subseteq \Gamma^{ij}(k|k)$
and thus
$\theta_{ij}(k) \in \Gamma^{ij}(k|k)$.
As $\theta_{ij}(0)\in \Gamma^{ij}(0|0)$, the thesis follows by induction.
\end{pf}
\vspace{-0.4cm}
\begin{lem}\label{lem:th_in_H_b}
Let Assumptions 1-4 hold. Then, for all $k$ such that $\mathcal{I}_i(k)=0$, $\theta_{i,i-1}(k) \in \Gamma^{ij}(k|k)$.
\end{lem}
\vspace{-0.7cm}
\begin{pf}
From \eqref{eq:u_est}, we have $\hat u^i_{i-1}(k) = [0, \ \omega_0+K]$ $\forall k:\mathcal{I}_i(k)=0$. Hence, from \eqref{eq:contr_strat} we conclude that $u_{i,i-1}(k)\in \hat u_{i,i-1}(k)$. The thesis then follows from Lemma \ref{lem:th_in_H}.  
\end{pf}
\vspace{-0.4cm}
\begin{defi}
We say that agent $i$ has reached the {\it desired spacing} with respect to agent $i-1$ at a generic time instant $k$ if $\mathcal{I}_i(k)=1$ and $\hat \vartheta_{i,i-1}(k)-\psi >0$ for the first time. 
\end{defi}
\vspace{-0.2cm}
Notice that, according to \eqref{eq:our_control_law_others_a}, when agent $i$ reaches the desired spacing, control is deactivated and $u_i(k)=\omega_0$.

\vspace{-0.4cm}
\section{Convergence Analysis}\label{sec:conv}
\vspace{-0.4cm}
Let us define the set
\vspace{-0.3cm}
\begin{equation}\label{eq:S}
\mathcal{S}:=\left\{k:\sum_{i=1}^N\mathcal{I}_i(k)\neq 0\right\}.
\end{equation}
\vspace{-0.3cm}
Now, we can state the following lemma.
\begin{lem}\label{thm:prel_2}
Let Assumptions 1-4 hold. If $\mathcal{S}\neq\emptyset$, then $\mathcal{I}_2(k_2)=1$, $\mathcal{I}_i(k_2)=0$ for all $i\ne 2$, where $k_2 = \min \mathcal{S}$.
\end{lem}
\vspace{-0.6cm}
\begin{pf}
If $u_{ij}(k) = 0$, then neither $i$ nor $j$ can discern if the other preceeds or follows as distance measurements give no information on orientation \cite{DeGa:15}. As $\mathcal{I}_i(0) = 0 \ \forall i$, we have that $u_{i,i-1}(0)= 0$ for all pairs of consecutive agents except the pairs $(2,1)$ and $(1,N)$, yielding $u_{21}(0)=-\omega_0$ and $u_{1N}(0)=\omega_0$, respectively. Hence, at $k=1$ the only agent that may discern its follower is agent $2$. This is still true for all $k$ such that $\mathcal{I}_i(k) = 0 \ \forall i$.
\end{pf}
\vspace{-0.6cm}
Now, we prove that agent $2$ identifies its follower in finite time.
\vspace{-0.2cm}
\begin{thm}\label{thm:idfol_ag2}
Let Assumptions 1-4 hold. Then, $\mathcal{S}\ne \emptyset$, and $k_2\leq\tilde{k}:= \left \lceil \frac{\theta_{21}(0)-2(\omega_0+K)}{\omega_0}\right \rceil$.
\end{thm}
\vspace{-0.6cm}
\begin{pf}
Notice that $\mathcal{S}\ne \emptyset$ is equivalent to the existence of $k_2$. Therefore we will prove the existence of $k_2$ and that it is smaller than $\tilde{k}$.
At $k=0$, $\mathcal{I}_i(k)=0$ $\forall i$.
From Lemma \ref{thm:prel_2}, we have that either $\mathcal S=\emptyset$, and therefore
$\mathcal{I}_i(k)=0$ for all $i$, or, if $\mathcal S\ne \emptyset$, then for all $k$
such that $\mathcal{I}_2(k)=0$, we also have $\mathcal{I}_i(k)=0$ for all $i\ne 2$. Summing up, we have that $\mathcal{I}_i(k)=0$ for all $i\ne 2$ until $\mathcal{I}_2(k)$ will become 1, if it ever happens. Now, to prove the thesis, it suffices to show that (i) \eqref{eq:follower_a} and (ii) \eqref{eq:follower_b} hold at time $\tilde{k}$.
\newline \noindent (i) From \eqref{eq:our_control_law_others_a} and \eqref{eq:our_control_law_leader}, we know that $u_{21}(k)=-\omega_0$ and $u_{2j}(k)=0$ for all $k$ such that $\mathcal{I}_2(k)=0$.
Let us generalize \eqref{eq:a_priori} as
\vspace{-0.3cm}
\begin{equation}\label{eq:a_priori_delta}
\Gamma^{ij}(k+\delta|k)=\Gamma_{ij}(k|k)+\sum_{\kappa=k}^{k+\delta-1}\hat{u}_{ij}(\kappa),
\end{equation}
where $\delta\in\mathbb{N}$. Now, observing that $\ushort\Gamma_l^{ij}(k|k)\ge\ushort\Gamma_l^{ij}(k|0)$, combining \eqref{eq:u_est} and \eqref{eq:a_priori_delta} we obtain
$ \ushort \Gamma_1^{2j}(k|k)\ge \ushort \Gamma_1^{2j}(k|0) \geq \ushort \Gamma_1^{2j}(0|0)\geq 2\varphi +2(\omega_0 + K)$
for all $j\ge 2$ and for all $k\leq \tilde k$.
From Lemmas \ref{lem:th_in_H} and \ref{lem:th_in_H_b}, eq. \eqref{eq:a_post_nonid_b}, and as $\vartheta_{21}(\tilde k)< 2(\omega_0+K)$, we have 
$\bar \Gamma_1^{21}(\tilde k|\tilde k) \leq \bar \Upsilon_{21}(\tilde k) < 2(\omega_0 + K)+2\varphi,$
and thus $\bar \Gamma_1^{21}(\tilde k|\tilde k)< \ushort \Gamma_1^{2j}(\tilde k|\tilde k) \ \forall j\ge 2$. Hence, \eqref{eq:follower_b} is fulfilled for $i=2$ at time $\tilde{k}$.
\newline \noindent (ii) From Lemma \ref{lem:th_in_H_b} and eqs. \eqref{eq:u_est}, \eqref{eq:our_control_law_others_a}, and \eqref{eq:our_control_law_leader}, $u_{21}(k) \in \hat u_{21}(k) = [0, \omega_0+K]$ for all $k \leq \tilde k$ and thus, as  $-\ushort \Upsilon_{21}(\tilde k) > -2(\omega_0 +K) - 2\varphi$ and  from \eqref{eq:a_priori} we have $\bar \Gamma_2^{21}(\tilde k|\tilde k) \leq -\theta_{21}(0) + 2\varphi < -\ushort \Upsilon_{21}(\tilde k)$ and thus 
$\Gamma_2^{21}(\tilde k|\tilde k)\cap (\Upsilon_{21}(\tilde k) \cup -\Upsilon_{21}(\tilde k))= \emptyset.$
Hence, $\ushort H^{21}(\tilde k|\tilde k) = \ushort \Gamma_2^{21}(\tilde k|\tilde k)>0$ and therefore
\eqref{eq:follower_a} holds at $\tilde k$ for $i=2$.
\end{pf}
\begin{table}
{\fontsize{8}{8}\selectfont
  \begin{center}
    \begin{tabular}{c|c}
    \scriptsize      \textbf{Symbols} & \scriptsize \textbf{Definition}\\
      \hline
      \scriptsize $k_i$ & \scriptsize Def. \ref{def:follower} \\
      \scriptsize $\mathcal{I}_i(k)$ & \scriptsize in eq. \eqref{eq:follower}\\
      \scriptsize $\hat{u}_i, \omega_0, K$  & \scriptsize in \eqref{eq:u_est}\\
      \scriptsize $\hat{u}_{ij}$ & \scriptsize after \eqref{eq:u_est}\\
      \scriptsize $\Gamma^{ij}$ & \scriptsize in \eqref{eq:est_strat}\\
      \scriptsize $\mathcal{S}$ & \scriptsize in \eqref{eq:S} 
    \end{tabular}
  \end{center}
  \rule{0pt}{3ex}\textsc{Table I.} Main symbols used in the all lemmas and theorems, and their definition in the main text.
  }
\end{table}
\begin{remark}\label{rema:th_2}
As $0<\vartheta_{21}(k_2)\leq\theta_{\max}$, $\psi> \theta_{\max}$, and from \eqref{eq:our_control_law_others_a}, 
we have that $u_{21}(k_2) = K$. Moreover, for all $k \geq k_2$, we have that the uncertainty set $\Gamma^{21}(k|k)$ is an interval. Namely, $\Gamma^{21}(k|k) = H^{21}(k|k) = \Gamma^{21}_1(k|k)$.
\end{remark}
\vspace{-0.2cm}
\begin{thm}\label{thm:err_bound}
Let Assumptions 1-4 be satisfied. Then, there exists a finite time $k_2^c$, such that $\vartheta_{21}(k)=\bar{\theta}_{21}$ for all $k\ge k_2^c$, where
\begin{align*}
k_2^c \leq & \left\lceil \frac{\theta_{21}(0)-2(\omega_0+K)}{\omega_0}\right\rceil + \left\lfloor \frac{\theta_{\max} - (\omega_0 + K)}{K} \right\rfloor\\
 &+1+\left\lceil\frac{\psi - (\theta_{\max}+K)}{K}\right\rceil;
\end{align*}
and $|\bar{\theta}_{21}-\psi |\leq K$.
\end{thm}
\vspace{-0.5cm}
\begin{pf}
As $0<\vartheta_{21}(k_2)\leq\theta_{\max}$, from \eqref{eq:a_post_fol} and \eqref{eq:scal_est} we have that $\hat \vartheta_{21}(k_2)<\theta_{\max}<\psi$, and thus \eqref{eq:our_control_law_others_a} and \eqref{eq:our_control_law_leader} imply $u_{21}(k) = K$, for all  $k>k_2$ such that $\alpha_{21}(k)<\theta_{\max}$. Hence, from  \eqref{eq:rel_ang_pos} there exists a time instant
\begin{equation*}
\tilde k_2= k_2+\left\lfloor \frac{\theta_{\max} - \vartheta_{21}(k_2)}{K} \right\rfloor+1
\end{equation*}
such that $\alpha_{21}(\tilde k_2) \in I^c$, and $\alpha_{21}(\tilde k_2-1) \leq \theta_{\max}$. Therefore, as we know from \eqref{eq:our_control_law_leader} that $u_1(k) = \omega_0$ for all $k$, then  \eqref{eq:u_est} implies $u_{21}(k)\in\hat u_{21}(k)$ for all $k$, and thus from Lemma \ref{lem:th_in_H} and eq. \eqref{eq:a_post_fol} we have
$\vartheta_{21}(\tilde k_2-1) \leq \bar H^{21}(\tilde k_2 -1|\tilde k_2-1)\leq \theta_{\max}.$
Then, applying \eqref{eq:a_priori} to $\Gamma^{21}(\tilde k_2-1|\tilde k_2-1) = H^{21}(\tilde k_2 -1|\tilde k_2-1)$, from \eqref{eq:u_est}, \eqref{eq:our_control_law_others_a}, and \eqref{eq:our_control_law_leader} we obtain
\vspace{-0.2cm}
\begin{equation}\label{eq:thH_k_2}
\vartheta_{21}(\tilde k_2) \leq \bar H^{21}(\tilde k_2 |\tilde k_2 -1)\leq \theta_{\max} + K.
\end{equation}
Finally, as $\alpha_{21}(\tilde k_2) \in I^c$, we have
$ \vartheta_{21}(\tilde k_2) > \theta_{\max}, $
and thus
\begin{equation}\label{eq:obvious}
0\leq \bar H^{21}(\tilde k_2 |\tilde k_2 -1)-\vartheta_{21}(\tilde k_2)<K.
\end{equation}
Hence, as from the estimation rule in \eqref{eq:scal_est} we have that $\hat\vartheta_{21}(\tilde k_2) = \bar{H}^{21}(\tilde k_2)$, the estimation error is bounded by $K$. Moreover, as $\psi \in (\theta_{\max}, \pi]$, \eqref{eq:u_est} and \eqref{eq:our_control_law_others_a} imply that, for all $n\in \mathbb{N}$ such that $\hat \vartheta_{21}(\tilde k_2+n-1)<\psi$,
\begin{equation}\label{eq:2728_a}
u_{21}(\tilde k_2 +n-1) = \hat u_{21}(\tilde k_2 +n-1) = K.
\end{equation}
This has two relevant consequences. Firstly, from \eqref{eq:a_priori} and \eqref{eq:a_post_fol}, for all $n\in \mathbb{N}$ such that $\hat \vartheta_{21}(\tilde k_2+n-1)<\psi$, we have
$H_{21}(\tilde k_2 +n|\tilde k_2+n) = H_{21}(\tilde k_2 +n|\tilde k_2 +n-1)
= H_{21}(\tilde k_2 +n-1|\tilde k_2+n-1)+K$,
and thus
\begin{equation}\label{eq:2728_c}
H_{21}(\tilde k_2 +n | \tilde k_2 +n) = H_{21}(\tilde k_2|\tilde k_2) +nK.
\end{equation}
Secondly,
\begin{equation}\label{eq:2728_b}
\vartheta_{21}(\tilde k_2 +n) = \vartheta_{21}(\tilde k_2) + nK.
\end{equation}
Hence, from \eqref{eq:obvious}, $\forall n\in \mathbb{N}_0:\hat \vartheta_{21}(\tilde k_2+n-1)<\psi$, we have
$
\bar H_{21}(\tilde k_2 +n|\tilde k_2 +n) -  \vartheta_{21}(\tilde k_2 +n) = \bar H_{21}(\tilde k_2|\tilde k_2) +nK - 
\vartheta_{21}(\tilde k_2) -nK
= \bar H_{21}(\tilde k_2|\tilde k_2) -  \vartheta_{21}(\tilde k_2),
$
which implies
\begin{equation}\label{eq:est_err_b}
0\leq \bar H_{21}(\tilde k_2+n|\tilde k_2+n) -  \vartheta_{21}(\tilde k_2+n)<K
\end{equation}
for all $n\in \mathbb{N}_0$ such that $\hat \vartheta_{21}(\tilde k_2+n-1)<\psi$.
Now, take
\begin{equation}\label{eq:n_2_def}
n_2 := \left\lceil\frac{\psi - \bar H_{21}(\tilde k_2|\tilde k_2)}{K}\right\rceil.
\end{equation}
By definition of the ceil function, we have that $n_2$ is such that $\hat \vartheta_{21}(\tilde k_2+n_2-1)<\psi$ and thus, from \eqref{eq:est_err_b}, 
$ 0\leq\bar H_{21}(\tilde k_2 +n_2|\tilde k_2+n_2) -  \vartheta_{21}(\tilde k_2 +n_2)< K$.
Now, as from the definition of $n_2$ in \eqref{eq:n_2_def} $ 0\leq \bar H_{21}(\tilde k_2 +n_2|\tilde k_2+n_2) -  \psi \leq K,$ 
we obtain that
$ |\vartheta_{21}(\tilde k_2 +n_2)-\psi|\leq K.$
Finally, as $\hat \vartheta_{21}(\tilde k_2+n_2)\geq \psi$, \eqref{eq:our_control_law_others_a} implies $u_2(\tilde k_2+n_2) = \omega_0$. Thus, from \eqref{eq:u_est} and \eqref{eq:our_control_law_leader} it follows that $\hat u_{21}(\tilde k_2+n_2) = u_{21}(\tilde k_2+n_2) = 0$. From \eqref{eq:a_post_fol} we also have 
$\hat u_{21}(k) = u_{21}(k) = 0$
for all the $k>\tilde k_2+n_2$, and thus $\vartheta_{21}(k)$ converges in finite time to a value $\bar\theta_{21}$ such that
$|\bar \theta_{21}-\psi|< K.$
Setting $k_2^c = \tilde k_2 + n_2$, the proof of existence of $k_2^c$ follows. Now, let us prove that
\begin{align*}
k_2^c \leq & \left\lceil \frac{\theta_{21}(0)-2(\omega_0+K)}{\omega_0}\right\rceil + \left\lfloor \frac{\theta_{\max} - (\omega_0 + K)}{K} \right\rfloor\\
&+1 +\left\lceil\frac{\psi - (\theta_{\max}+K)}{K}\right\rceil.
\end{align*}
To do so, let us start by considering that $k_2^c = \tilde k_2 + n_2$, which from  \eqref{eq:n_2_def} implies that 
$k_2^c = \tilde k_2 + \left\lceil\frac{\psi - \bar H_{21}(\tilde k_2|\tilde k_2)}{K}\right\rceil.$
In turn, from \eqref{eq:obvious} and from the definition of $\tilde k_2$ we have
\begin{align*}
k_2^c &= \tilde k_2 + \left\lceil\frac{\psi - \theta_{\max}+K}{K}\right\rceil\\
 &=k_2+\left\lfloor \frac{\theta_{\max} - \vartheta_{21}(k_2)}{K} \right\rfloor+1+ \left\lceil\frac{\psi - \theta_{\max}+K}{K}\right\rceil.
\end{align*}
Finally, from Theorem 1 which ensures that $k_2\leq  \left \lceil \frac{\theta_{21}(0)-2(\omega_0+K)}{\omega_0}\right \rceil$
we obtain 
$
k_2^c \leq \left\lceil \frac{\theta_{21}(0)-2(\omega_0+K)}{\omega_0}\right\rceil + \left\lfloor \frac{\theta_{\max} - (\omega_0 + K)}{K} \right\rfloor
+1+\left\lceil\frac{\psi - (\theta_{\max}+K)}{K}\right\rceil.
$
\end{pf}
\vspace{-0.6cm}
Theorems \ref{thm:idfol_ag2} and \ref{thm:err_bound} prove that our strategy allows to bound the steady state value of $|\theta_{21}(k)-\psi |$ with $K$ which, being a parameter of the control law, can be made arbitrarily small either by slowing down the agents or by reducing the sampling time. Now, we extend such results to the remainder of the multi-agent system. To do so we will consider a generic agent $i$, and make some assumptions on agents $i-1$ and $i-2$. The following two Lemmas will prove that these assumptions guarantee the convergence of the estimation and control strategy converge, respectively. Then, we will prove that they are always verified for $i=3,...,N$.

\begin{lem}\label{thm:conv_stima}
Let Assumptions 1-4 hold. For all  $i=3,\ldots,N$, if $k_{i-1}$ exists, $\mathcal{I}_i(k_{i-1}) = 0$, and $\vartheta_{i-1,j}(k_{i-1})>\omega_0+K$, for all $ j \neq i : \vartheta_{i-1,j}(k_{i-1})>0$, then $k_i$ exists and $0\leq \vartheta_{i,i-1}(k_i) \leq \theta_{\max}$. 
\end{lem}
\vspace{-0.6cm}
\begin{pf}
To prove the thesis, it suffices to show that there exists a time instant $k_i$ when \eqref{eq:follower} holds. The proof is organized in two steps, where we show that \eqref{eq:follower_a} and \eqref{eq:follower_b} hold, respectively.
\newline \noindent \textit{Step 1.} To prove \eqref{eq:follower_a}, we distinguish between three cases:
\begin{enumerate}[1.]
\item $\vartheta_{i,i-1}(0) \in (\theta_{\max},\pi)$;
\item $4\varphi +2 \omega_0 + 2K < \vartheta_{i,i-1}(0) \leq \theta_{\max}$;
\item $\vartheta_{i,i-1}(0) < 0$.
\end{enumerate}
Note that in all cases, as from \eqref{eq:our_control_law_others_a} we have that $u_{i,i-1}(k) = 0 \ \forall k<k_{i-1}$, then $\vartheta_{i,i-1}(k_i) = \vartheta_{i,i-1}(0)$.
\newline Case 1. From \eqref{eq:a_post_nonid_a}, \eqref{eq:our_control_law_leader}, and \eqref{eq:our_control_law_others_c}, and from Lemma \ref{thm:prel_2}, we know that $u_{i,i-1}(k)<0$ for all $k\geq k_{i-1}$ such that $\vartheta_{i,i-1}(k)\in (\theta_{\max}, \pi)$. Hence, at time 
\begin{align*}
\tilde{k}&\leq k_{i-1} + \left \lceil \left(\vartheta_{i,i-1}(k_{i-1})-\theta_{\max}\right)/\omega_0\right\rceil\\
&= k_{i-1} + \left \lceil \left(\vartheta_{i,i-1}(0)-\theta_{\max}\right)/\omega_0\right\rceil
\end{align*}
we have that $\alpha_{i,i-1}(\tilde{k})\le \theta_{\max}$, while $\alpha_{i,i-1}(k)> \theta_{\max}$ for all $k<\tilde k$. From \eqref{eq:a_priori}-\eqref{eq:a_post_nonid_a}, we have that $\Gamma_2^{i,i-1}(\tilde{k}|\tilde{k}-1)=[-\pi,-\theta_{\max})$ and $-\ushort{\Upsilon}^{i,i-1}(\tilde k)\ge\theta_{\max}$ yielding  
$\Gamma_2^{i,i-1}(\tilde k|\tilde k-1) \cap -\Upsilon^{i,i-1}(\tilde k) = \emptyset $ and thus 
$\ushort{H}_{i,i-1}(\tilde k|\tilde k) = \ushort{\Gamma}_1^{i,i-1}(\tilde k|\tilde k)>0$. 
\newline Case 2.
In this case,
$
\vartheta_{i,i-1}(k_{i-1}) \in \Gamma^{i,i-1}(k_{i-1}|k_{i-1})
 = \Gamma_1^{i,i-1}(k_{i-1}|k_{i-1}) \cup \Gamma_2^{i,i-1}(k_{i-1})
$
where, as $\vartheta_{i,i-1}(k_{i-1}) =\vartheta_{i,i-1}(0)$, and thus $ \Gamma^{i,i-1}(k_{i-1}|k_{i-1}) \in \Gamma^{i,i-1}(0|0),$ we have that
\begin{align}
&\bar{\Gamma}_1^{i,i-1}(k_{i-1}|k_{i-1}) \leq \theta_{\max},\\
&\ushort{\Gamma}_1^{i,i-1}(k_{i-1}|k_{i-1}) \geq 2\varphi+\omega_0+K,\\
&\ushort{\Gamma}_2^{i,i-1}(k_{i-1}|k_{i-1}) \geq -\theta_{\max},\\
&\bar{\Gamma}_2^{i,i-1}(k_{i-1}|k_{i-1}) \leq -2\varphi-\omega_0-K.
\end{align}
Following the line of argument in Case 1, we could show that at a time 
$\tilde{k}\leq k_{i-1} + \left \lceil \frac{\vartheta_{i,i-1}(k_{i-1})-4\varphi}{\omega_0}\right\rceil$
we would have that $0<\vartheta_{i,i-1}(\tilde k)<\vartheta_{i,i-1}(k_i)-4\varphi$. From \eqref{eq:a_priori_delta}, we obtain $\bar{\Gamma}_2^{i,i-1}(k|0) = \bar{\Gamma}_2^{i,i-1}(0|0)$,
for all $k\le\tilde{k}$. Hence, 
\begin{equation}
\begin{aligned}
\bar{\Gamma}_2^{i,i-1}(\tilde k|0) \le -\vartheta_{i,i-1}(0)+2\varphi-\ushort{\Upsilon}^{i,i-1}(\tilde k) \ge\\ \ge -\vartheta_{i,i-1}(\tilde k)-2\varphi >-\vartheta_{i,i-1}(0)+2\varphi
\end{aligned}
\end{equation}
and therefore $\Upsilon^{i,i-1}(\tilde k) \cap \Gamma_2^{i,i-1}(\tilde k|0) = \emptyset$. Thus, we have 
$\ushort{H}_{i,i-1}(\tilde k|\tilde k) = \ushort{\Gamma}_1^{i,i-1}(\tilde k|\tilde k)>0.$
\newline Case 3. The proof can be conducted following the same steps as in Case 1, but setting
$\tilde{k}\leq k_{i-1}+ \left \lceil \frac{\theta_{i,i-1}(0)-\theta_{\max}}{\omega_0}\right\rceil.$
\newline \noindent \textit{Step 2.} Now, we prove that there exists a time instant in which \eqref{eq:follower_b} holds. Again, let us distinguish between:
\begin{itemize}
\item[] Case 1. $\vartheta_{i,i-1}(k_{i-1}) \in (\theta_{\max},\pi)$;
\item[] Case 2. $4\varphi +2 \omega_0 + 2K < \vartheta_{i,i-1}(k_{i-1}) \leq \theta_{\max}$;
\item[] Case 3. $\vartheta_{i,i-1}(k_{i-1}) < 0$.
\end{itemize} 
Case 1. If $\vartheta_{i,i-1}(k_{i-1}) \in (\theta_{\max},\pi)$, then \eqref{eq:a_post_nonid_a} and \eqref{eq:our_control_law_others_c} imply $u_i(k) = 0$ for all $k$ such that $\vartheta_{i,i-1}(k) \in (\theta_{\max},\pi)$. Indeed, from \eqref{eq:our_control_law_others_a}, there exists a time instant $\tilde k$ in which we will have that $\theta_{\max} - \omega_0 - K<\vartheta_{i,i-1}(\tilde k)\leq \theta_{\max}$, $\vartheta_{i,i-1}(\tilde k-1)>\theta_{\max}$ where
$\tilde k$ is defined as in Case 1 of Step 1. By hypothesis, and as at time $\tilde k$ we have that $\theta_{\max} - \omega_0 - K<\vartheta_{i,i-1}(\tilde k)$, then $\nexists j\neq i-1:0<\vartheta_{ij}(\tilde k)\leq \theta_{\max},$ and thus \eqref{eq:follower_b} holds for all $j$ such that $0<\vartheta_{ij}(\tilde k)$. Hence, consider the case in which, at time $\tilde k$, \eqref{eq:follower_b} has not been verified yet for an agent $j$ such that $\vartheta_{ij}(\tilde k)<0$. Indeed, as from \eqref{eq:our_control_law_others_c} we have that $u_i(k)=0 \ \forall k< \tilde k$, then from Assumption 2 we have
$\ushort \Gamma^{ij}_1(\tilde k|\tilde{k})\geq 2\varphi+2\omega_0+2K.$
Hence, \eqref{eq:follower_b} will hold before a time instant $\tilde k + n_i$ such that $\vartheta_{i,i-1}(\tilde k+n_i)\in [\omega_0+K, 2\omega_0+2K)$ where 
$n_i\leq \left\lfloor \frac{\theta_{\max} - (2\omega_0 +2K)}{\omega_0}\right\rfloor+1.$
In this case, as $u_i(k)=0 \ \forall k<\tilde k + n_i$ we would have that also $u_{ij}(k)\leq0 \ \forall k<\tilde k + n_i$. Hence, from Assumption 2, eqs. \eqref{eq:a_post_nonid_b} and \eqref{eq:c_iniz},  and as the width of $\Upsilon$ is bounded by $2\varphi$, we obtain that
$\bar H^{i,i-1}(\tilde k +n_i) < \ushort \Gamma^{ij}_1(\tilde k),$
thus verifying \eqref{eq:follower_b}.
\newline Case 2. From Assumption 2, eqs. \eqref{eq:a_post_nonid_b} and \eqref{eq:c_iniz}, and as the width of $\Upsilon$ is bounded by $2\varphi$, we have that \eqref{eq:follower_b} holds at time $k=0$ for all $j$ such that $\vartheta_{ij}(0)\in [0,\pi)\cup[-\pi, -\theta_{\max})$. The proof that \eqref{eq:follower_b} will eventually hold also for all agents such that $\vartheta_{ij}(0)\in[-\theta_{\max},0)$ can be performed following the same arguments made above and setting $\tilde k = k_{i-1}$.
\newline Case 3. The proof can be completed as in case 1, but noting that 
$\tilde{k}\leq k_{i-1}+ \left \lceil \left(2\pi + \vartheta_{i,i-1}(0)-\theta_{\max}\right)/\omega_0\right\rceil.$
\end{pf}
\begin{lem}\label{thm:conv_contr}
Let Assumptions 1-4 hold. For all $i=3,\ldots,N$, If there exists $h_i \geq k_i$ such that $0<\vartheta_{i,i-1}(h_i)<\theta_{\max}$, $\vartheta_{i,i-1}(h_i)\in \Gamma^{i,i-1}(h_i|h_i)$, and $\psi - \hat \vartheta_{j,j-1}\leq 0$ for all $j = 2,...,i-1$, then there exists $k_i^c:|\vartheta_{i,i-1}(k_i^c)-\psi |\leq K$ for all $k \geq k_i^c$.
\end{lem}
\vspace{-0.4cm}
\begin{pf}
As $\psi - \hat \vartheta_{i-1,i-2}\leq 0$, from \eqref{eq:our_control_law_others_a} we have that $u_{i-1}(k)=\omega_0$ $\forall k \geq h_i$. Moreover, as $h_i\geq k_i$ and $0<\vartheta_{i,i-1}(h_i)<\theta_{\max}$, \eqref{eq:a_post_fol} and \eqref{eq:scal_est} imply that $\hat \vartheta_{i,i-1}(h_i)<\theta_{\max}<\psi$, and thus from  \eqref{eq:our_control_law_others_a} we have $u_{i,i-1}(k) = K$ for all  $k>h_i$ such that $\alpha_{i,i-1}(k)<\theta_{\max}$. Hence, from \eqref{eq:rel_ang_pos} there exists a time instant
$\tilde k_i= h_i + \left\lfloor \left(\theta_{\max} - \vartheta_{i,i-1}(h_i)\right)/K \right\rfloor+1$ 
such that $\alpha_{i,i-1}(\tilde k_i) \in I^c$, and $\alpha_{i,i-1}(\tilde k_i-1) \leq \theta_{\max}$.  Therefore, as $u_{i-1}=\omega_0 \ \forall k\geq h_i$, by hypothesis and from \eqref{eq:u_est}, $u_{i,i-1}(k)\in \hat u_{i,i-1} \ \forall k$, and thus, from Lemma \ref{lem:th_in_H} and \eqref{eq:a_post_fol}, we have $\vartheta_{i,i-1}(\tilde k_i-1)\leq \bar H^{i,i-1}(\tilde k_i -1|\tilde k_i -1)\leq \theta_{\max}$. Then, applying \eqref{eq:a_priori} to $\Gamma^{i,i-1}(\tilde k_i-1|\tilde k_i -1) = H^{i,i-1}(\tilde k_i -1|\tilde k_i -1)$, from \eqref{eq:u_est} and as $\alpha_{i,i-1}(\tilde k_i)\in I^c$, we have 
$\theta_{\max}<\vartheta_{i,i-1}(\tilde k_i) \leq \bar H^{i,i-1}(\tilde k_i |\tilde k_i -1)\leq \theta_{\max} + K,$
which, thanks to the estimation rule in \eqref{eq:scal_est}, ensures the estimation error is bounded by $w(H^{i,i-1}(\tilde k_i |\tilde k_i))< K$. Moreover, as $\psi \in (\theta_{\max}, \pi]$, from \eqref{eq:u_est} and \eqref{eq:our_control_law_others_a}, for all $n\in \mathbb{N}_+$ such that $\hat \vartheta_{i,i-1}(\tilde k_i+n-1)<\psi$ we have
$u_{i,i-1}(\tilde k_i +n-1) = \hat u_{i,i-1}(\tilde k_i +n-1) = K,$
which, in turn, has two relevant implications. Firstly, from \eqref{eq:a_post_fol} and \eqref{eq:our_control_law_others_a}, it implies that $H^{i,i-1}(\tilde k_i +n|\tilde k_i+n) = H^{i,i-1}(\tilde k_i +n-1|\tilde k_i+n-1)+K$
and thus, from 
\eqref{eq:a_post_fol} and \eqref{eq:a_priori_delta},
\begin{equation}\label{eq:h_kn_k}
H^{i,i-1}(\tilde k_i +n|\tilde k_i+n) = H^{i,i-1}(\tilde k_i|\tilde k_i)+nK.
\end{equation}
Secondly, it implies that
\begin{equation}\label{eq:vartheta_kn}
\vartheta_{i,i-1}(\tilde k_i +n) = \vartheta_{i,i-1}(\tilde k_i) + nK.
\end{equation}
Subtracting \eqref{eq:vartheta_kn} from \eqref{eq:h_kn_k}, we have that $\forall n\in \mathbb{N}$ such that $\hat \vartheta_{i,i-1}(\tilde k_i+n-1)<\psi$, $
\bar H^{i,i-1}(\tilde k_i+n|\tilde k_i+n) -\vartheta_{i,i-1}(\tilde k_i+n) =\bar H^{i,i-1}(\tilde k_i|\tilde k_i) -\vartheta_{i,i-1}(\tilde k_i)$
and from \eqref{eq:scal_est}
\begin{align}
|\hat \vartheta_{i,i-1}(\tilde k_i +n) - \vartheta_{i,i-1}(\tilde k_i+n)| = \nonumber\\\label{eq:err_fin}
=|\hat \vartheta_{i,i-1}(\tilde k_i) - \vartheta_{i,i-1}(\tilde k_i)|< K.
\end{align}
Now, take the time instant
$n_i := \left\lceil\frac{\psi - \bar H^{i,i-1}(\tilde k_i|\tilde k_i)}{K}\right\rceil.$
From the definition of the ceil function, we have that
\begin{equation}\label{eq:h_theta}
0\leq \bar H^{i,i-1}(\tilde k_i+n_i|\tilde k_i+n_i)-\vartheta_{i,i-1}(\tilde k_i+n_i)<K,
\end{equation}
while from the definition of $n_i$ we obtain that
\begin{equation}\label{eq:h_psi}
0\leq \bar H^{i,i-1}(\tilde k_i+n_i|\tilde k_i+n_i)-\psi\leq K.
\end{equation}
Combining \eqref{eq:h_theta} and \eqref{eq:h_psi}, we obtain
$|\vartheta_{i,i-1}(\tilde k_i+n_i) - \psi|\leq K.$
Finally, at time $\tilde k_i+n_i$, as we have that $\hat \vartheta_{i,i-1}(\tilde k_i+n_i) = \bar H^{i,i-1}(\tilde k_i+n_i|\tilde k_i+n_i)\geq \psi$, then from  \eqref{eq:u_est} and \eqref{eq:our_control_law_others_a} we have that $\hat u_{i,i-1}(\tilde k_i +n_i) = u_{i,i-1}(\tilde k_i +n_i) = 0$. From  \eqref{eq:a_post_fol}, this is also true for all $k\leq \tilde k_i +n_i$, and thus $\vartheta_{i,i-1}(k)$ converges in finite time to a value $\bar \theta_{i,i-1}$ such that $ |\bar \theta_{i,i-1}-\psi|<K \ \forall k\geq \tilde k_i +n_i.$ Setting $k_i^c = \tilde k_i +n_i$, the thesis follows.
\end{pf}
\vspace{-0.3cm}
Lemmas \ref{thm:conv_stima} and \ref{thm:conv_contr} prove convergence of both the estimation and control strategies under some given assumptions. Hence, we now only need to prove that the hypotheses therein are always verified for each agent $i$.
\begin{thm}\label{thm:final_conv}
If Assumptions 1-4 hold and $K<\epsilon/(N-1)$, then the proposed estimation and control strategy is capable of achieving an $\varepsilon$-balanced circular formation. Moreover, for all $(i,j)\in\left\{ (1,2),\ldots,(N-1,N)\right\},$ 
\begin{equation}\label{eq:conv_ij}
\lim_{k\rightarrow\infty}\left|\theta_{ij}(k)-\psi\right|\le K.
\end{equation}
\end{thm}
\vspace{-0.8cm}
\begin{pf}
From Theorem \ref{thm:err_bound}, we know that $|\vartheta_{21}(k)-\psi|\leq K$ $\forall k\geq k_2^c$. To prove the thesis for $i=3$, we must first prove the hypotheses of Lemmas \ref{thm:conv_stima} and \ref{thm:conv_contr} hold for $i=3$. Let us start from Lemma \ref{thm:conv_stima}, that is by proving that $k_{2}$ exists, $\mathcal{I}_3(k_2) = 0$, and $\vartheta_{3,j}(k_2)>\omega_0+K$, for all $ j \neq 3 : \vartheta_{2,j}(k_{2})>0$.
\newline \noindent The existence of $k_2$ is proven in Theorem \ref{thm:idfol_ag2} and from Lemma \ref{thm:prel_2} we know that $\mathcal{I}_3(k_2) = 0$. Moreover, as $\mathcal{I}_i(k_2) = 0$ for all $i\neq 2$, then, from \eqref{eq:our_control_law_others_a} and Assumption 2, $|\vartheta_{3j}(k_2)| = |\vartheta_{3j}(0)|\geq 4\varphi+2\omega_0+2K$ for all $j\neq 1$.  Finally, from Remark \ref{rema:th_2}, we also have that $\vartheta_{21}(k_2)>0$ and thus $\vartheta_{31}(k_2)>\vartheta_{32}(k_2)$ and thus the hypotheses of Lemma \ref{thm:conv_stima} hold for $i=3$ which implies $k_3$ exists. 
\newline \noindent Now, let us prove that the hypothesis of Lemma \ref{thm:conv_contr} hold for $i=3$, that is, that the time instant $h_3$ exists. Define the time instant $\tilde h_3 := \max \lbrace k_3, k^c_2\rbrace$. From Lemma \ref{thm:conv_stima}, we know that $0\leq \vartheta_{i,i-1}(k_3)\leq \theta_{\max}$ and that as $\mathcal{I}_{2}(k_3)=1$, if $k_2^c>k_3$, then from eq. \eqref{eq:our_control_law_others_a} we have $\vartheta_{32}(k) = \vartheta_{32}(k_3)$ for all $k$ in $[k_3 \ k_2^c]$. Hence, $0\leq \vartheta_{i,i-1}(\tilde h_3)\leq \theta_{\max}$. Moreover, from Lemma \ref{lem:th_in_H_b} we know that $\theta_{32}(k_3) \in \Gamma_{32}(k_3)$ and if $k_3<k^2_c$, from eq. \eqref{eq:our_control_law_others_c} $u_{32}(k) = 0$, and as we have proved that $\vartheta_{32}(k) = \vartheta_{32}(k_3)$ for all $k$ in $[k_3 \ k_2^c]$, from eq. \eqref{eq:u_est} we have $u_{32}(k)\in \hat u_{32}(k)$. Finally, from Lemma \ref{lem:th_in_H}, we have that $\theta_{32}(\tilde h_3)\in \Gamma_{32}(\tilde h_3)$, and, as $\psi - \hat \vartheta_{21}(\tilde h_3)\leq 0$, the definition of $\tilde h_3$ implies $h_3 =\tilde h_3$. Hence, the existence of $k_2$, which is guaranteed by Theorem \ref{thm:idfol_ag2} ensures the existence of $k_3$, while the existence $k_2^c$, guaranteed by Theorem \ref{thm:err_bound} ensures the existence of $k_3^c$. If we prove that $k_4>k_{3}$ then this reasoning could be iterated for all pairs $i$ and $i+1$ starting from the pair $(3,4)$ to the pair $(N-1,N)$ and thus \eqref{eq:conv_ij} would follow by induction. Following the same line of arguments of Lemma \ref{thm:conv_stima}, it is possible to prove that if $0<\vartheta_{3,2}\leq 2(\omega_0+k)$, then $\mathcal{I}_3(k)=1$. From \eqref{eq:contr_strat}, this ensures that for all $k$ such that $\vartheta_{43}(k)$ and $\vartheta_{42}(k)$ are both greater than zero we have that $\vartheta_{42}(k)>\vartheta_{43}(k)$. In turn, from Lemma \ref{lem:th_in_H_b} this implies that $\bar H^{42}(k)>\ushort \Gamma^{43}_1(k)$ for all $k$, and thus, from \eqref{eq:follower} the closest follower of agent $4$ can only be agent $3$. Hence, considering that \eqref{eq:our_control_law_others_c} $\vartheta_{43}(k) = \vartheta_{43}(0)$ for all $k$ such that $\mathcal{I}_3(k)=0$ and that from \eqref{eq:u_est} and \eqref{eq:a_priori}-\eqref{eq:a_post_nonid_b}  both $\Gamma^{43}_1(k)$ and $\Gamma^{43}_2(k)$ are nonempty, we have that \eqref{eq:follower_a} cannot hold for $i=4$ before $k_3$. Hence, \eqref{eq:conv_ij} follows by induction. Now, as agent $1$ travels at constant speed $\omega_0$ and thus, for all $k\geq k_{N}^c$, we have that
$\vartheta_{1,N}(k)-\psi = 2\pi -\sum_{i=2}^N \vartheta_{i,i-1} \in
 [2\pi-(N-1)(\psi+K), 2\pi-(N-1)(\psi-K)]$
and, as $N\psi = 2\pi$, then $\vartheta_{1,N}(k) = [\psi -(N-1)K,\psi +(N-1)K],$
which finally implies $\epsilon$-bounded convergence, with $\epsilon=(N-1)K$.
\end{pf}
The following theorem completes the results of Theorem \ref{thm:err_bound} by providing upper bounds for the convergence times $k_i^c$ of the estimation and control strategy for the remaining agents $i=3,\ldots,N$.
\begin{thm}\label{thm:t_conv_bound}
Let Assumptions 1-4 hold. For all $i = 3,\dots,N$, if $\vartheta_{i,i-1}(0)>0$, then
\begin{align*}
k_i^c \leq &\max \left\lbrace k_{i,i-1} + \left\lceil \frac{\vartheta_{i,i-1}(0) - 2(\omega_0 + K)}{\omega_0}\right\rceil, k_{i-1}^c\right\rbrace\\
&+ \left\lceil\left(\psi - 2(\omega_0+K)\right)/K\right\rceil.
\end{align*}
Otherwise,
\begin{align*}
k_i^c \leq &\max \left\lbrace k_{i,i-1} + \left\lceil \frac{(2\pi + \vartheta_{i,i-1}(0)) - 2(\omega_0 + K)}{\omega_0}\right\rceil, \right.\\
& k_{i-1}^c\bigg\rbrace+ \left\lceil\left(\psi - 2(\omega_0+K)\right)/K\right\rceil.
\end{align*}
\end{thm}
\begin{pf}
Iterating the reasoning performed in Theorem \ref{thm:final_conv}, we can prove that $k_i>k_{i-1}$. Hence, $\vartheta_{i,i-1}(k_{i-1}) = \vartheta_{i,i-1}(0)$. Then, as $u_{i,i-1}(k) \leq \omega_0$ for all $k$ such that $k_{i-1}\leq k < k_i$, and as following the lines of argument of theorem \ref{thm:idfol_ag2} for agent $2$ we can prove that $\omega_0+K\leq \vartheta_{i,i-1}(k_i)\leq \theta_{\max}$, if $\vartheta_{i,i-1}(0)>0$, we have
\begin{equation}\label{eq:bound1}
k_i \leq k_{i,i-1} + \left\lceil \left(\vartheta_{i,i-1}(0) - 2(\omega_0 + K)\right)/\omega_0\right\rceil.
\end{equation}
Otherwise, we have that 
\begin{equation}\label{eq:bound2}
k_i \leq k_{i,i-1} + \left\lceil \left((2\pi + \vartheta_{i,i-1}(0)) - 2(\omega_0 + K)\right)/\omega_0\right\rceil.
\end{equation}
Then, as (i) the first time instant in which $u_{i,i-1}(k)>0$ is $\max \lbrace k_i, k_{i-1}^c\rbrace,$ and if $k_{i-1}^c>k_{i}$, then $\vartheta_{i,i-1}(k_{i-1}^c) = \vartheta_{i,i-1}(k_{i})$;
(ii) from Theorem \ref{thm:err_bound} $|\vartheta_{i,i-1}(k_i^c) - \psi|<K$, then for all $k$ such that $\max \lbrace k_i, k_{i-1}^c\rbrace\leq k<k_i^c$ we have $u_{i,i-1}(k) = K$. 
Combining (i) and (ii) yields $
k_i^c = \max \lbrace k_i, k_{i-1}^c\rbrace + \left\lceil \left(\psi - \vartheta_{i,i-1}(\max \lbrace k_i, k_{i-1}^c\rbrace)\right)/K\right\rceil
\leq \max \lbrace k_i, k_{i-1}^c\rbrace + \left\lceil\left(\psi - 2(\omega_0+K)\right)/K\right\rceil.$
Substituting $k_i$ with one of the two bounds derived in \eqref{eq:bound1} and \eqref{eq:bound2}, the thesis follows.
\end{pf}
\begin{figure}
\centering
\includegraphics[scale=0.40, trim = 1cm 8cm 0cm 8cm]{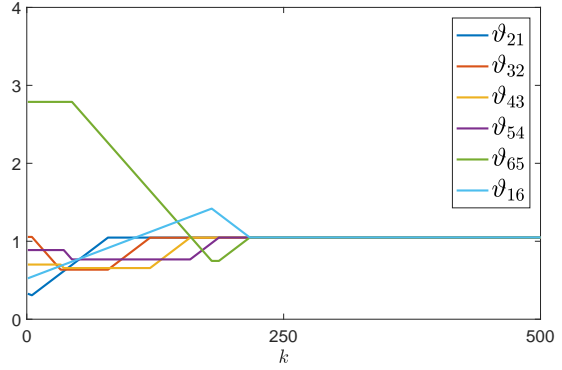}
\caption{Time evolution of the phase differences between consecutive agents in a sample simulation.}\label{fig:example}
\end{figure}
\begin{table}
{\fontsize{8}{8}\selectfont
\rule{0pt}{1ex}  
\begin{center}
\begin{tabular}{c|c|c|c|c}
 & $\mathbf{K = \omega_0}$ & $\mathbf{K = 2\omega_0}$ & $\mathbf{K = 3\omega_0}$ & $\mathbf{K = 4\omega_0}$ \\ 
\hline 
$\mathbf{\langle k^c_6 \rangle}$ & 748 & 434 & 339 & 287 \\ 
\hline 
$\mathbf{\eta}$ &  3.1$\scriptstyle\times$10$^{-3}$ & $6.4\scriptstyle\times$10$^{-3}$ & $7.3\scriptstyle\times$10$^{-3}$ & $14.2\scriptstyle\times$10$^{-3}$
\end{tabular}
\end{center}
\rule{0pt}{4ex}\textsc{Table II.} Variation of the average convergence time $\langle k^c_6 \rangle$ and of the average steady state error $\eta= \langle |\bar \theta_{i,i-1} - \psi|\rangle$ as a function of the control gain $K$.
}
\end{table}
\section{Numerical validation}\label{sec:num}
We consider $N=6$ agents, which implies $\psi = \pi/3$. Moreover, we set the value of $\theta_{\max}$ to $3/4\psi= \pi/4$, $\omega_0$ to $0.005$, and let the control gain $K$ take the values in the set $\lbrace \omega_0, \ 2\omega_0, \ 3\omega_0, \ 4\omega_0\rbrace$. Finally, for each value of $K$, we vary $\varphi$ in the set $\lbrace 2K, \ 3K, \ 4K, \ 5K\rbrace$, and select, for the random variable $\nu_{ij}(k)$, a uniform distribution in the interval $[-\varphi \ \varphi]$. Such parameter selection defines 16 different scenarios for each of which we run $100$ numerical experiments where the initial conditions are taken randomly in the admissible region of the state-space defined by Assumption 2. Figure \ref{fig:example} shows the plot of the time evolution of $\vartheta_{i,i-1}(k)$ for all $i$ for a representative simulation. The numerical results are consistent with our theoretical predictions, as in all cases the multi-agent system achieved an $\varepsilon$-balanced circular formation with $\varepsilon = K(N-1)$ consistently with Theorem \ref{thm:final_conv}. In coherence with the latter, we observe that $|\bar \theta_{i,i-1}-\psi|<K$ for all $i = 3,\dots,6$. Consistently with the bound derived in Theorem \ref{thm:t_conv_bound}, we observe a trade-off between the convergence time $k^c_6$ and the accuracy of the balanced formation $\epsilon$ that depends on the value of $K$. As the gain $K$ increases, the speed of convergence also increases
($\langle k^c_6 \rangle$ reduces), but the average steady state error $\eta:=\langle |\bar \theta_{i,i-1} - \psi| \rangle $ increases, see Table II.
\section{Conclusions}\label{sec:concl}
We proposed an estimation and control strategy for balancing a formation of autonomous agents on a circle in the case in which only proximity measurements with a radius that is lower than the desired spacing are available, implying that the agents are blind when approaching the desired formation. This setting reproduces situations in which only inexpensive proximity sensors can be employed, and few agents can be deployed to patrol a given boundary. We exploit the limited information coming from the measurement equation through a three-level bang bang controller that is symbiotic with our estimation strategy. Our completely decentralized approach prescribes a random election of a pacemaker, which sets the pace of the system. We showed that the system achieves an $\varepsilon$-partially bounded circular formation, and that this bound can be made arbitrarily small by leveraging a control parameter. The theoretical analysis is complemented by a set of simulations illustrating the trade-off between the speed of convergence and the accuracy of the formation balancing. Future work will be devoted to extend our approach to cope with different agent dynamics.

\bibliographystyle{plain}

\end{document}